\documentclass[12pt]{amsart}
\usepackage{amscd,amssymb}
\newtheorem{theorem}{Theorem}[section]
\newtheorem{lemma}[theorem]{Lemma}
\newtheorem{claim}[theorem]{Claim}
\newtheorem{proposition}[theorem]{Proposition}
\newtheorem{corollary}[theorem]{Corollary}

\theoremstyle{definition}
\newtheorem{definition}[theorem]{Definition}

\theoremstyle{remark}
\newtheorem{remark}[theorem]{Remark}

\theoremstyle{conj}
\newtheorem{conjecture}[theorem]{Conjecture}

\begin{document}

\today

\title[non-singular solution]
{Non-singular solutions to the normalized Ricci flow equation}

\author[F. Fang]{Fuquan Fang}
\thanks{The first author was supported by
NSF Grant 19925104 of China, 973 project of Foundation Science of
China, and the Capital Normal University}
\address{Nankai Institute of Mathematics,
Weijin Road 94, Tianjin 300071, P.R.China}
\address{Department of Mathematics, Capital Normal University,
Beijing, P.R.China}
  \email{ffang@nankai.edu.cn}
\author[Y. Zhang]{Yuguang Zhang}
\address{Department of Mathematics, Capital Normal University,
Beijing, P.R.China  }
\author[Z. Zhang]{Zhenlei Zhang}
\address{Nankai Institute of Mathematics,
Weijin Road 94, Tianjin 300071, P.R.China}

\begin{abstract}

In this paper we study non-singular solutions of Ricci flow on a
closed manifold of dimension at least $4$. Amongst others we prove
that, if $M$ is a closed $4$-manifold on which the normalized
Ricci flow exists for all time $t>0$ with uniformly bounded
sectional curvature, then the Euler characteristic $\chi (M)\ge
0$. Moreover, the $4$-manifold satisfies one of the following

\noindent (i) $M$ is a shrinking Ricci solition;

\noindent (ii) $M$ admits a positive rank $F$-structure;

\noindent (iii) the Hitchin-Thorpe type inequality holds
$$ 2\chi (M)\ge 3|\tau(M)|$$
where $\chi (M)$ (resp. $\tau(M)$) is the Euler characteristic
(resp. signature) of $M$.
\end{abstract}
\maketitle



\section{Introduction}

In this paper we will consider the normalized Ricci flow equation
on a closed smooth $n$-dimensional manifold $M$,
\begin{equation}\label{0}
\left\{ \begin{array}{ll}
\frac{\partial}{\partial t}g=-2Ric+\frac{2r}{n}g, \\
g(0)=g_{0},
\end{array} \right.
\end{equation}
where $Ric$ and $r$ denote the Ricci tensor and average scalar
curvature $\frac{\int_{M}Rdv}{\int_{M}dv}$ of the Riemannian
metric $g$ respectively. Following Hamilton \cite{H}, a solution
to this equation is called {\it non-singular}, if the flow exists
for all time $t\geq0$ and the Riemannian curvature tensor
satisfies $|Rm|\leq C<\infty$ uniformly for some constant $C$
independent of $t$.

In \cite{H}, Hamilton classified $3$-dimensional non-singular
solutions. In particular, he proved that the underlying 3-manifold
is geometrizable in the sense of Thurston. Hamilton's theorem is
of great importance to understand the long-time behavior of
solutions to the Ricci flow, even to the Ricci flow with surgery
which was used by Perelman in \cite{P1} as a technique to study
the global property of the Ricci flow solution, modulo modifying
the singular points in space time (cf. for example \cite{CZ} or
\cite{KL} for detailed discussion.)

Our main result is the following generalization of Hamilton's
results to higher dimensions.
\begin{theorem}\label{001}
Any non-singular solution to the normalized Ricci flow equation
(\ref{0}) on a closed manifold $M$ does one and only one of the
following things:
\begin{itemize}\label{00}
 \item[(1.1.1)] the solution collapses;
 \item[(1.1.2)] the solution converges along a subsequence to a
 shrinking Ricci soliton solution on $M$;
 \item[(1.1.3)] the solution converges along a subsequence to a
 Ricci flat metric solution on $M$;
 \item[(1.1.4)] the solution converges along a subsequence to a
 negative Einstein metric solution on $M$;
 \item[(1.1.5)] the solution converges along a subsequence of times and space points to a complete non-compact
negative Einstein metric solution in the pointed Gromov-Hausdoff
sense.
\end{itemize}
\end{theorem}

A solution $(M,g(t))$ to equation (\ref{0}) is called collapse if
the maximum of the local injectivity radius
$\text{inj}(x,g(t)),x\in M,$ tends to zero along some subsequence
of times $t_{k}\rightarrow\infty$.

In dimension 3, by a result of Ivey \cite{I}, every closed Ricci
soliton is Einstein. So the convergence result coincides with that
of Hamilton's. For case (1.1.5), Hamilton proved that the
$3$-manifold splits into pieces of hyperbolic manifolds and
residual graph manifolds. Each hyperbolic piece has finite volume
and finite cusps at infinity, every cusp contains an
incompressible torus $T$, $i.e.$, $\pi_{1}(T)$ injects into
$\pi_{1}(M)$, with constant mean curvature and small area. In
general, by Cheeger-Gromov [3], the Riemannian manifold $(M,g(t))$
for large $t$ admits a thick-thin decomposition
$M=M^{\varepsilon}\bigcup M_{\varepsilon}$ for small
$\varepsilon$, where
$$M^{\varepsilon}=\{x\in M|\text{Vol}(B(x,1,g(t)))\geq\varepsilon\},$$
$$M_{\varepsilon}=\{x\in M|\text{Vol}(B(x,1,g(t)))\leq\varepsilon\}.$$
By (1.1.5) it is easy to see that the thick part is recognized as
the negative Einstein pieces, while the thin part has an
$F$-structure which partially collapses when we take the limit
(cf. \cite[Thm. 1.3]{A}). In 4-dimension we will prove the thin
part is indeed volume collapsed (cf. Theorem 1.5).

The following is a generalization of the celebrated Hitchin-Thorpe
inequality for non-singular solutions to the Ricci flow on closed
$4$-manifolds.

 \vskip 3mm

\begin{theorem} \label{002}
Let $M$ be a closed oriented  4-manifold $M$, and
$\{g(t)\},t\in[0,\infty)$, be a non-singular  solution to
(\ref{0}).
 Then $M$   satisfies one of the following

\noindent (1.2.1) $M$ admits a shrinking Ricci solition;

\noindent (1.2.2) $M$ admits a positive rank $F$-structure;

\noindent (1.2.3) the Hitchin-Thorpe type inequality holds
$$ 2\chi (M)\ge 3|\tau(M)|$$
where $\chi (M)$ (resp. $\tau(M)$) is the Euler characteristic
(resp. signature) of $M$.
\end{theorem}

By [8] a closed shrinking Ricci soliton has finite fundamental
group. Thus, a $4$-dimensional closed shrinking Ricci soliton has
Euler characteristic at least $2$. On the other hand, manifold
with positive rank $F$-structure has vanishing Euler
characteristic. Therefore Theorem 1.2 implies readily that

\vskip 3mm

\begin{corollary}\label{003}
If a  closed 4-manifold $M$ admits a non-singular solution, then
$\chi(M)\geq0$.
\end{corollary}

The above corollary gives a topological obstruction to the
existence of non-singular solutions on a closed $4$-manifold.

To state a strengthened result of Theorem 1.2 we need to introduce
Perelman's $\lambda$-functional (cf. \cite{P1}\cite{KL}\cite{FZ}).

For a smooth function $f\in C^\infty(M)$ on a Riemannian
$n$-manifold with a Riemannian metric $g$, let
\begin{equation}\mathcal{F}(g,f)=\int_{M} (R_g+|\nabla f|^{2})e^{-f}
dvol_{g},\end{equation} where $R_g$ is the scalar curvature of
$g$.
 The Perelman's
$\lambda$-functional  is defined by
\begin{equation}\lambda_{M}(g)=\inf _f \{\mathcal{F}(g,f)|\int_{M}e^{-f}
 dvol_{g}=1\}.\end{equation} Note that $\lambda_{M}(g)$ is the lowest eigenvalue of the
 operator $-4\triangle+R_{g}$.  Let \begin{equation}\overline{\lambda}_{M}(g)
 =\lambda_{M}(g)\text{Vol}_{g}(M)^{\frac{2}{n}}\end{equation}
which is invariant up to rescale the metric.  Perelman \cite{P1}
has established the monotonicity property of $\overline{\lambda}
_M(g_t)$ along the Ricci flow $g_t$, namely, the function is
non-decreasing along the Ricci flow $g_t$ whenever
$\overline{\lambda} _M(g_t)\leq 0$. Therefore, it is interesting
to study the upper bound of $\overline{\lambda}_{M}(g)$. This
leads to define a diffeomorphism invariant
  $\overline{\lambda}_{M} $ of $M$ due to Perelman (cf. \cite{P2}
  \cite{KL})
by \begin{equation}\overline{\lambda}_{M}=
  \sup\limits_{g\in \mathcal{M}} \overline{\lambda}_{M}(g),
  \end{equation} where $\mathcal{M}$ is the set of Riemannian
metrics on $M$. As we pointed out in \cite{FZ} $\overline{\lambda}
_M$ is not a homeomorphism invariant. Observe that
$\overline{\lambda} _M=0$ if $M$ admits a volume collapsing with
bounded scalar curvature but does not admit any metric with
positive scalar curvature (cf. \cite{KL}).

A upper bound in \cite{FZ} for the invariant
$\overline{\lambda}_M$ was obtained by using the Seiberg-Witten
theory, whenever $\overline{\lambda}_M< 0$ and $M$ has a
non-trivial monopole class.

Now we are ready to state

\begin{theorem}\label{004}
Let $M$ be a closed oriented  4-manifold $M$ with
$\overline{\lambda}_{M}< 0$, and $\{g(t)\},t\in[0,\infty)$, be a
solution to (\ref{0}).  If $|R(g(t))|<C$ where $C$ is a constant
independent of $t$, then
$$2\chi(M)-3|\tau(M)|\geq \frac{1}{96\pi^{2}}\overline{\lambda}_{M}^{2},$$
In particular, if $\{g(t)\},t\in[0,\infty)$ is a non-singular
solution to (\ref{0}), then the solution  $(M,g(t))$ does not
collapse.
\end{theorem}

The above theorem combined with the Seiberg-Witten theory implies
that
 \vskip 3mm

\begin{corollary} \label{005}
Let $M$ be a closed symplectic 4-manifold $M$ with
$\overline{\lambda}_{M} <  0$, and $\{g(t)\},t\in[0,\infty)$, be a
solution to (\ref{0}) with $|R(g(t))|<C$ where $C$ is a constant
independent of $t$.
 If $b^{+}_{2}(M):=\text{dim}H^2_+(M;\Bbb R)
   >1$, then $$\chi(M) \geq 3\tau(M).$$
\end{corollary}

In dimension $4$ Theorem (1.1.5) may be improved as follows:

\vskip 3mm

\begin{theorem} \label{006}
 Let $(M,g(t)),t\in[0,\infty)$ be a non-singular
solution to (\ref{0}) on a closed oriented  4-manifold $M$ with
$\overline{\lambda}_{M} < 0$. Then, for any $\delta >0$, there is a
time $T\gg 1$, and a compact 4-submanifold $M^{\varepsilon}$ with
boundary in $M$, $M^{\varepsilon}\subset M$, such that
\begin{itemize}\label{00}
 \item[(1.5.1)]  $\text{Vol}(M-M^{\varepsilon}, g(T))< \delta$, and $M-M^{\varepsilon}$ admits an F-structure of positive rank.
  \item[(1.5.2)] The
 components of $\partial M^{\varepsilon}$ are graph 3-manifolds. \item[(1.5.3)] $M^{\varepsilon}$ admits
 an Einstein metric with negative scalar curvature $g_{\infty}$ which is close to
 $g(T)|_{M^{\varepsilon}}$ in the $C^{\infty}$-sense.
  \end{itemize}
\end{theorem}

In view of Theorem 1.6, it is natural to wonder what kind of
information can be obtained on the thick part $M^\varepsilon$ in
the above theorem. Applying the Seiberg-Witten theory we obtain
the following result which provides a partial answer for certain
symplectic manifolds.

\vskip 3mm

\begin{theorem}\label{006} Let $(M,g(t))$ and  $(M^{\varepsilon},g_{\infty})$ be the same  as in Theorem 1.6.  If $M$
  admits a symplectic structure satisfying that $b^{+}_{2}(M)
   >1$ and  $\chi(M)=3 \tau (M)$, then $g_{\infty}$ is
   a complex hyperbolic metric.
\end{theorem}

\vskip 3mm

We conclude this section by posing the following

\begin{conjecture}
Theorem (1.2.3) may be replaced by the following
Hitchin-Thorpe-Gromov-Kotschik type inequality
$$2\chi (M)-3|\tau (M)|\ge \frac 1{1296\pi ^2}\|M\|$$
 where $\|M\|$ is a simplicial volume of $M$.
\end{conjecture}

 The organization of the paper is as follows: In $\S$2 we give a proof of
Theorem 1.1. In $\S$3 we are concerned with $4$-dimensional
non-singular solutions, and we will prove Theorem 1.2, Theorem 1.4
and Corollary 1.5. In $\S 4$ we will prove Theorem 1.6. Finally,
we will prove Theorem 1.7 in $\S 5$.

\vskip 10mm


\section{Proof of Theorem 1.1}

In \cite{H2}, Hamilton introduced the (unnormalized) Ricci flow
equation
\begin{equation}\label{101}
\left\{ \begin{array}{ll} \frac{\partial}{\partial t}g=-2Ric, \\
g(0)=g_{0},
\end{array} \right.
\end{equation}
and its normalized equation (\ref{0}). By \cite{H2} the equation
(\ref{0}) is just a change of equation (\ref{101}) via rescalings
in space and a reparametrization in time, such that the volumes of
the Riemannian metrics are preserved to be constant. In this
paper, we will always assume the volumes of the metrics equal $1$,
whenever we consider the solutions to the normalized Ricci flow.
Especially we have $\text{Vol}(g_{0})=1$.

As in \cite{P1}, one also can define a scale invariant version of
$\mathcal{F}$ functional, the so called $\mathcal{W}$ functional.
For each smooth $f$ and constant $\tau>0$, let
\begin{equation}
\mathcal{W}(g,f,\tau)=\int_{M}[\tau(R+|\nabla
f|^{2})+f-n](4\pi\tau)^{-n/2}e^{-f}dv,\nonumber
\end{equation}
and then set
\begin{equation}
\mu(g,\tau)=\inf\{\mathcal{W}(g,f,\tau)|\int_{M}(4\pi\tau)^{-n/2}e^{-f}dv=1\},\nonumber
\end{equation}
\begin{equation}
\nu(g)=\inf\limits_{\tau>0}\mu(g,\tau).\nonumber
\end{equation}
 By a result
of Rothaus \cite{R}, for each $\tau>0$, there exists a smooth
minimizer $\phi$ such that $\mu(g,\tau)=\mathcal{W}(g,\phi,\tau)$.
By Claim 3.1 of \cite{P1}, if $\lambda(g)>0$, then the infimum in
the definition of $\nu$ is attained by some $\tau>0$ and
$\nu(g)\leq0$. By \cite{P1} the $\lambda$ functional and the $\nu$
functional are non-decreasing along the Ricci flow, which plays a
central role in our proof (cf. Lemma \ref{105} and Lemma \ref{106}
for proofs of these facts.)

We start with the definition of convergence of solutions to the
equation (\ref{0}) or (\ref{101}). We assume that all the
Riemannian metrics taking into account are complete.
\begin{definition}[\cite{H1}]
Let $(M_{k},g_{k}(t),p_{k}), p_{k}\in M, t\in(A,B)$ with
$-\infty\leq A<0$ and $0<B\leq\infty$, be a sequence of marked
solutions to the Ricci flow equation (\ref{0}) or (\ref{101}). We
say that $(M_{k},g_{k}(t),p_{k})$  converges to another solution
$(M_{\infty},g_{\infty}(t),p_{\infty}),t\in(A,B),$ to (\ref{0}) or
(\ref{101}) respectively, if there is a sequence of increasing
open subsets $U_{k}\subset M_{\infty}$ containing $p_{\infty}$,
$i.e.$, $U_{k}\subset U_{k+1}$ for each $k$, and a sequence of
diffeomorphisms $F_{k}:U_{k}\to V_{k}\subset M_{k}$ mapping
$p_{\infty}$ to $p_{k}$, such that the pull-backed metrics
$F_{k}^{*}g_{k}(t)$ converge to $g_{\infty}(t)$ on every compact
subset of $M_{\infty}\times(A,B)$ uniformly together with all
their derivatives.
\end{definition}

In \cite{H1}, Hamilton proved his famous compactness theorem for
solutions to the unnormalized Ricci flow equation (\ref{101}),
under the assumption: (i) the local injectivity radii of the
metrics $g_{k}(0)$ at $p_{k}$
 are uniformly bounded below;
(ii) the supremum norm of Riemannian curvature tensors are
uniformly bounded above on any compact time interval. For a local
version of this theorem, see \cite{CZ} for example. We remark that
the compactness theorem of Hamilton remains valid for the
normalized Ricci flow solutions under the same assumption. This
can be checked easily from the proof of the theorem.

Recall that $(M,g(t))$ is a {\it Ricci solition} solution to the
Ricci flow, if $g(t)$ is obtained from $g(0)$ via changes by
diffeomorphisms and rescalings. At each time the solution
satisfies $$Ric+\mathcal{L}_{X}g=cg$$ for some vector field $X$
and some constant $c$. The Ricci soliton is said to be shrinking,
steady or expanding according to $c>0,c=0,c<0$ respectively. Note
that if the manifold is closed, then $c=r$, which can be seen by
taking the trace of above equation and then integrating it over
$M$, where $r$ is the average scalar curvature of $g$.

In \cite{H}, Hamilton considered metrics with positive scalar
curvature. Our next proposition deals with the solutions so that
the $\lambda$-functional is positive.

\vskip 3mm

\begin{proposition}\label{102}
Let $(M,g(t)),t\in[0,\infty),$ be a non-singular solution to
(\ref{0}) on a closed $n$-manifold $M$ with $\lambda(g(0))>0$.
Then for any sequence of times $t_{k}\rightarrow\infty$, there is
a subsequence $t_{k_{i}}$ such that $(M,g(t_{k_{i}}+t))$ converges
to a shrinking Ricci soliton solution.
\end{proposition}
\begin{proof}
If we denote by $(M,\bar{g}(t))$ the corresponding unnormalized
Ricci flow solution, then by Proposition 1.2 of \cite{P1},
$$\frac{d}{dt}\lambda(\bar{g}(t))\geq\frac{2}{n}\lambda(\bar{g}(t))^{2}.$$
Thus the unnormalized Ricci flow must extinct in finite time and
$\lambda(\bar{g}(t))>0$ remains hold. Since $g(t)$ are just
rescalings of $\bar{g}(t)$, $\lambda(g(t))>0$ is preserved. Then
Perelman's no local collapsing theorem \cite[Chap. 4]{P1} tells us
that there exist $\kappa,\rho>0$ such that each metric ball
$B=B_{g(t)}(p,r)\subset M$, with radius $r\leq\rho$ and
$\sup_{x\in B}|Rm|(x,t)\leq r^{-2}$, has volume
$\text{Vol}_{g(t)}(B)\geq\kappa r^{n}$. By the definition of
non-singular solution, there is a constant $C>0$ such that
$|Rm|(x,t)\leq C$ for all $(x,t)\in M\times[0,\infty)$. By
replacing $\rho$ by a smaller constant, we may assume $\rho\leq
C^{-\frac{1}{2}}$. Thus $|Rm|(x,t)\leq\rho^{-2}$ always hold and
so $Vol(B_{g(t)}(x,\rho))\geq\kappa\rho^{n}$ for all $(x,t)$. This
implies that the diameters of $g(t)$ are bounded above uniformly,
since we can fill only $\frac{1}{\kappa\rho^{n}}$ disjoint balls
of radius $\rho$ in $M$ at each time $t$. By a result of Cheeger,
Gromov and Taylor \cite{CGT}, the injectivity radii of $g(t)$ are
bounded below uniformly for any $t$. Then Hamilton's compactness
theorem yields the convergence result.

Let $(M,g_{k}(t)=g(t_{k}+t))$, $t_{k}\rightarrow\infty$, be such a
sequence which converges to a limit solution to the normalized
Ricci flow equation
$(M_{\infty},g_{\infty}(t)),t\in(-\infty,\infty)$. By the uniform
boundedness of diameters of $g(t)$, $M_{\infty}=M$ and the
convergence is smooth on any compact time interval. Next we will
show that $(M,g_{\infty}(t))$ is a shrinking Ricci soliton.

 Now for each
time $t\in(-\infty,\infty)$, $\nu(g(t))$ is achieved by some
positive number $\tau(t)$, $i.e.$,
$\nu(g(t))=\mu(g(t),\tau(t))\leq0$ holds. Further, by Lemma
\ref{106} below,  $\nu(g(t))$ increases. By the smooth convergence
of $(M,g_{k}(t))$ for any $t$, we have
$$\nu(g_{\infty}(t))=\lim\limits_{k\rightarrow\infty}\nu(g(t_{k}+t))
=\lim\limits_{k\rightarrow\infty}\nu(t_{k}+t)
=\lim\limits_{t\rightarrow\infty}\nu(t),$$ is a constant
independent of $t$. We consider two
subcases for the limit solution. \\
Case 1: $\lambda(g_{\infty}(t_{0}))>0$ for some $t_{0}\geq0$. In
this case, $\nu(g_{\infty}(t))$ is attainable for $t$ around
$t_{0}$ and using Lemma \ref{106} again, one sees that
$(M,g_{\infty}(t))$ is really a Ricci soliton. By definition of
$\lambda$, the average scalar curvature
$r(g_{\infty}(t_{0}))\geq\lambda(g_{\infty}(t_{0}))>0$, so the
Ricci
soliton is a shrinking one.\\
Case 2: $\lambda(g_{\infty}(t))\equiv0$. In this case, the
monotonicity of the $\lambda$ functional implies that
$(M,g_{\infty}(t))$ is Einstein. See Lemma \ref{105} below. The
scalar curvature is zero because it equals
$\lambda(g_{\infty}(t))$. We will exclude this exception by
showing that $\nu(g_{\infty}(t))=-\infty$, which is a
contradiction since
$\nu(g_{\infty}(t))=\lim_{t\rightarrow\infty}\nu(g(t))\geq\nu(g(0))$
by the increasing of $\nu(g(t))$ along the Ricci flow. Hence the
proposition from the following claim.
\end{proof}

\begin{claim}
Let $(N,g)$ be a closed Riemannian manifold with average scalar
curvature $r\leq0$. Then we have
$\lim\limits_{\tau\rightarrow\infty}\mu(g,\tau)=-\infty$.
\end{claim}
\begin{proof}[Proof of the claim]
For any smooth function $f$, set $u=(4\pi\tau)^{-n/2}e^{-f}$; then
by the definition of the $\mathcal{W}$ functional,
$$\mathcal{W}(g,f,\tau)=\int_{M}[\tau(Ru+\frac{|\nabla u|^{2}}{u})-u\ln u]dv-\frac{n}{2}\ln(4\pi\tau)-n.$$
Choosing $u=1$ and substituting it into the functional, we have
\begin{eqnarray}
\mu(g,\tau)+n&\leq&\tau
r-\frac{n}{2}\ln(4\pi\tau)\leq-\frac{n}{2}\ln(4\pi\tau)\rightarrow-\infty\nonumber
\end{eqnarray}
as $\tau\rightarrow\infty$. This ends the proof of the claim.
\end{proof}

From the proof of above proposition, we have an immediate
corollary.
\begin{corollary}
If $(M,g)$ is a closed Riemannian manifold with $\lambda(g)>0$,
then the normalized Ricci flow solution, with $g$ as initial
metric, will never converge to a Ricci flat metric.
\end{corollary}

By a result of Fern\'{a}ndez-L\'{o}pez and Garc\'{i}a-R\'{i}o
\cite{FG}, any shrinking Ricci soliton has finite fundamental
group, so we have
\begin{corollary}
Let $(M,g)$ be a closed Riemannian manifold with $\lambda(g)>0$.
If the normalized Ricci flow solution on $M$ with $g$ as initial
metric is non-singular, then $\pi_{1}(M)$ is finite.
\end{corollary}

For a Riemannian metric $g$ on a manifold $M$, denote by
$\mbox{inj}(x,g)$ the injectivity radius of the metric $g$ at $x$.
We say a solution to the normalized Ricci flow collapses if there
exists a sequence of times $t_{k}\rightarrow T$ such that
$\sup_{x\in M}\mbox{inj}(x,g(t_{k}))\rightarrow0$, where $T$ is
the maximal existence time for the solution, which may be finite
or infinite. If a non-singular solution to the normalized Ricc
flow equation doesn't collapse, then by Hamilton's compactness
theorem, for each sequence of times $t_{k}\rightarrow\infty$,
there exists a subsequence $t_{k_{i}}$ and a sequence of points
$p_{i}\in M$ such that $(M,g(t_{k_{i}}+t),p_{i})$ converges to
another normalized Ricci flow solution
$(M_{\infty},g_{\infty}(t),p_{\infty})$. Denote by
$r_{\infty}(t)=\lim_{i\rightarrow\infty}r(g(t_{k_{i}}+t))$ the
limit constant in the normalized Ricci flow equation, then we have
$$\frac{\partial}{\partial t}g_{\infty}(t)=-2Ric_{\infty}(t)+\frac{2}{n}r_{\infty}(t)g_{\infty}(t),$$
where $Ric_{\infty}(t)$ is the Ricci tensor of $g_{\infty}(t)$.
Note that $r_{\infty}(t)$ may not equal to the average scalar
curvature $r(g_{\infty}(t))$.

Denote by $\breve{R}(g)=\min_{x\in M}R(x)$ the minimum of the
scalar curvature of a given metric $g$.

\vskip 2mm

\begin{proposition}\label{103}
Let $(M,g(t)),t\in[0,\infty)$, be a non-singular solution to
(\ref{0}) on a closed $n$-manifold $M$. Assume that
$\breve{R}(g(t))\leq-c<0$ uniformly for some constant $c>0$
independent of $t$. If the solution doesn't collapse, then there
exists a sequence of points $p_{k}\in M$ and a sequence of times
$t_{k}\rightarrow\infty$ such that the solutions
$\{(M,g(t_{k}+t),p_{k})\}_{k=1}^{\infty}$ converge to an Einstein
metric solution to (\ref{0}), whose scalar curvature is negative.
\end{proposition}

For proving this proposition, we need the following lemma:

\begin{lemma}\label{107} Let $(M,g(t)),t\in[0,\infty)$ be a  solution to
(\ref{0}) on a closed $n$-manifold $M$. Assume that
$\breve{R}(g(t))\leq-c<0$ uniformly for some constant $c>0$
independent of $t$. If $|R(g(t))|\leq C$, where $C$ is a constant
independent of $t$, then
$$\int_{0}^{\infty}(r(g(t))-\breve{R}(g(t)))dt<\infty, \ \ \ {\rm and} \ \
\
$$ $$\int_{0}^{\infty}\int_{M}|R(g(t))-r(g(t))|dvdt<\infty.$$
\end{lemma}

\begin{proof}
We follow the proof given in Section 7 of \cite{H} by Hamilton.
Consider the evolution equation of $R$
\begin{eqnarray}
\frac{\partial}{\partial t}R=\triangle
R+2|Ric\textordmasculine|^{2}+\frac{2}{n}R(R-r),\nonumber
\end{eqnarray}
where $Ric\textordmasculine$ denotes the traceless part of Ricci
tensor. By maximal principle,
$\frac{d}{dt}\breve{R}\geq\frac{2}{n}\breve{R}(\breve{R}-r)$ and
so $\breve{R}(g(t))$ increases whenever it is negative. By
assumption
$\frac{d}{dt}\breve{R}\geq\frac{2}{n}\breve{R}(\breve{R}-r)\geq\frac{2c}{n}(r-\breve{R})$,
which implies that
$$\int_{0}^{\infty}(r-\breve{R})dt<\infty.$$

Come back to the original solution $(M,g(t))$. We have
$$\int_{M}|R-r|dv\leq\int_{M}(R-\breve{R})dv+\int_{M}(r-\breve{R})dv=2\int_{M}(r-\breve{R})dv=2(r-\breve{R}).$$
Thus
$$\int_{0}^{\infty}\int_{M}|R-r|dvdt<\infty.$$
\end{proof}

\begin{proof}[Proof of Proposition \ref{103}]
Assume that
$\lim_{t\rightarrow\infty}\breve{R}(t)=-\delta\leq-c<0$. Note that
by assumption, there exists a constant $C>0$ such that
$|Rm|(x,t)\leq C$ uniformly for all $(x,t)\in M\times[0,\infty)$.
First we will show that $\lim_{t\rightarrow\infty}r(t)$ exists and
equals to $-\delta$. For this, consider the evolution equation
$$\frac{d}{dt}(r-\breve{R})\leq\int_{M}(2|Ric|^{2}-\frac{2}{n}Rr)dv-\frac{2}{n}\breve{R}(\breve{R}-r)\leq D,$$
where $D\geq1$ is a constant depending only on $C$ and $n$. We
claim that for any $0<\epsilon<1$, there is $T$ such that
$r(t)-\breve{R}(t)<\epsilon$ whenever $t>T$. Otherwise, there will
be a sequence of times $t_{k}\rightarrow\infty$ satisfying
$t_{k+1}\geq t_{k}+1$ and $r(t_{k})-\breve{R}(t_{k})\geq\epsilon$.
From the above equation, we obtain that
$r(t)-\breve{R}(t)\geq\frac{\epsilon}{2}$ whenever
$t\in[t_{k}-\epsilon/D,t_{k}]$. Thus
$$\int_{0}^{\infty}(r-\breve{R})dt\geq\sum_{k=1}^{\infty}\int_{t_{k}-\epsilon/D}^{t_{k}}(r-\breve{R}(t))dt
\geq\sum_{k=1}^{\infty}\int_{t_{k}-\epsilon/D}^{t_{k}}\frac{\epsilon}{2}dt=\infty,$$
which contradicts the Lemma \ref{107}. Hence
$\lim_{t\rightarrow\infty}r(t)=-\delta$ and consequently
$r_{\infty}(t)=-\delta$ on any limit solution as we mentioned in
the paragraphs before Proposition \ref{103}.

Considering the time interval $[t_{0},t_{0}+1],\forall
t_{0}\in{\mathbb{R}},$ on the limit solution, we have
$\int_{t_{0}}^{t_{0}+1}\int_{M_{\infty}}|R(g_{\infty}(t))-r_{\infty}(t)|dvdt=0$.
Hence $R(g_{\infty}(t))\equiv r_{\infty}(t)=-\delta$. Then by the
evolution of scalar curvature we have
$Ric\textordmasculine(g_{\infty}(t))=0$, $i.e.$,
$(M_{\infty},g_{\infty}(t))$ is Einstein for all $t$. The scalar
curvature
$R_{\infty}(t)=\lim_{t\rightarrow\infty}\breve{R}(t)=-\delta<0$.
\end{proof}

\begin{remark}
It is remarkable that the conclusion in Proposition \ref{103}
remains valid if we replace the boundedness of the Riemannian
curvature tensor by the boundedness of the Ricci tensor. This can
be seen from the process of proving the Proposition \ref{103}.
Also note that in both of these special cases, the limit constant
$r_{\infty}(t)$ equal to $r(g_{\infty}(t))$, the average scalar
curvature of the limit metrics.
\end{remark}

By now for a solution to the normalized Ricci flow equation, if
$\lambda(g(t))>0$ for some time $t$, then we can use Proposition
\ref{102}; while if $\breve{R}(g(t))\leq-c<0$ holds for all $t$,
then we can use Proposition \ref{103}. As for the remaining case,
we have the following proposition similar as the zero sectional
curvature limit case considered by Hamilton \cite{H}.
\begin{proposition}\label{104}
Let $(M,g(t)),t\in[0,\infty)$, be a non-singular solution to
(\ref{0}) on a closed $n$-manifold $M$. Assume that
$\breve{R}(g(t))\leq0$, $\lambda(g(t))\leq0$ and
$\breve{R}(g(t))\nearrow0$ as $t\rightarrow\infty$. If the
solution doesn't collapse, then there exists a sequence of times
$t_{k}\rightarrow\infty$ such that the solutions
$\{(M,g(t_{k}+t))\}_{k=1}^{\infty}$ converge to a Ricci flat
metric solution to (\ref{0}).
\end{proposition}
\begin{proof}
We divide the proof into several subcases.

Case 1: The unnormalized Ricci flow extincts in finite time.
Perelman's no local collapsing theorem shows that the limit
manifold is closed and so $M_{\infty}=M$. Now
$\lambda(g(t))\nearrow0$ implies that $\lambda(g_{\infty}(t))=0$
for all $t$. By Lemma \ref{105}, $(M,g_{\infty}(t))$ is a Ricci
flat solution.

Case 2: The unnormalized Ricci flow exists for all time
$t\in[0,\infty)$ and there is a sequence of times
$t_{k}\rightarrow\infty$ such that the average scalar curvature
$r(t_{k})$ of $g(t_{k})$ converges to zero. In this case, we use a
modified version of the proof used in Section 6 of \cite{H} by
Hamilton. Denote by $R_{\infty}$ and $Ric_{\infty}$ the scalar
curvature and Ricci tensor of the limit solution respectively. By
assumption,
\begin{equation}
\int_{M}(R(t_{k})-\breve{R}(t_{k}))dv_{g(t_{k})}=r(t_{k})-\breve{R}(t_{k})\rightarrow0,\nonumber
\end{equation}
as $k\rightarrow\infty$. Note that
$R(t_{k})-\breve{R}(t_{k})\geq0$. By taking the limit, one obtains
$\int_{M_{\infty}}R_{\infty}(0)dv_{g_{\infty}}(0)=0$. But
$R_{\infty}(0)\geq0$ over $M_{\infty}$ since
$\lim_{k\rightarrow\infty}\breve{R}(t_{k})=0$, so
$R_{\infty}(0)\equiv0$. Note that $r_{\infty}(0)=0$, since
$r(t_{k})\rightarrow0$. Consider the evolving equation of the
scalar curvature on the limit solution
\begin{equation}\nonumber
\frac{\partial}{\partial t}R_{\infty}(t)=\triangle
R_{\infty}(t)+2|Ric_{\infty}(t)|^{2}-\frac{2}{n}r_{\infty}(t)R_{\infty}(t),t\in(-\infty,\infty).
\end{equation}
It follows from the strong maximal principle that
$R_{\infty}\equiv0$ and $Ric_{\infty}\equiv0$ over
$M_{\infty}\times(-\infty,\infty)$, $i.e.$,
$(M_{\infty},g_{\infty}(t))$ is a Ricci flat solution. Now the
volume of the limit manifold is less than or equals to 1. By a
result of Yau \cite{SY}, $M_{\infty}$ is compact. So
$M_{\infty}=M$ and the convergence is smooth.

Case 3: The unnormalized Ricci flow, say $(M,\bar{g}(t))$, exists
for all time $t\in[0,\infty)$ and the average scalar curvature of
normalized Ricci flow $r(t)\geq\delta$ uniformly for some constant
$\delta>0$. We want to show this case will never happen. Denote by
$\overline{V}(t)=\text{Vol}(\bar{g}(t))$ the volume of
$\bar{g}(t)$. Since $r(g)\text{Vol}(g)^{\frac{2}{n}}$ is scale
invariant, we have
$\bar{r}(t)\geq\delta\overline{V}^{\frac{-2}{n}}(t)$. By the
evolving equation
\begin{equation}\nonumber
\frac{d}{dt}\overline{V}=\int_{M}-R(\bar{g})dv_{\bar{g}}=-\bar{r}\overline{V}\leq-\delta\overline{V}^{\frac{n-2}{n}},
\end{equation}
we obtain $\overline{V}^{\frac{2}{n}}(t)\leq1-\frac{2}{n}\delta t$
for all $t$, which contradicts with the assumption that the
unnormalized Ricci flow solution exists for all time. The desired
result follows.
\end{proof}

\begin{remark}
In fact, by a refined argument using the monotonicity of
Perelman's $\mu$ functional along the Ricci flow, Case 1 of
Proposition \ref{104} also can be excluded. So the only phenomena
is that of Case 2, under the assumption in the proposition.
\end{remark}

Summing up the results of Proposition \ref{102}, \ref{103} and
\ref{104}, we finish the proof of Theorem \ref{001}.

To conclude this section, let us prove two lemmas used previously,
which are basically due to Perelman. These are basic facts in the
study of Ricci flow solutions.

\vskip 3mm

\begin{lemma}[Perelman]\label{105}
Let $(M,g(t)),t\in[0,T),$ be a solution to the normalized Ricci
flow equation (\ref{0}) on a closed $n$-manifold $M$. Then
$\lambda(g(t))$ increases whenever $\lambda(g(t))\leq0$, and the
increasing is strict unless $g(t)$ is Einstein.
\end{lemma}
\begin{proof}
Consider the coupled equation
\begin{equation}
\left\{ \begin{array}{ll}
\frac{\partial}{\partial t}g=-2Ric+\frac{2r}{n}g, \\
\frac{\partial}{\partial t}f=-\triangle f+|\nabla f|^{2}-R+r.
\end{array} \right.\nonumber
\end{equation}
Under this evolving equation, we have
\begin{eqnarray}
\frac{d}{dt}\mathcal{F}(g(t),f(t))&=&\int_{M}[2|Ric+\nabla^{2}f|^{2}-\frac{2}{n}r(R+\triangle
f)]e^{-f}dv.\nonumber
\end{eqnarray}
Let $\bar{f}$ be the eigenfunction of $-4\triangle_{g(t)}+R(g(t))$
. Denote by $\lambda(t)=\lambda(g(t))$, then
\begin{eqnarray}
\frac{d}{dt}\lambda(t)&=&\int_{M}[2|Ric+\nabla^{2}\bar{f}|^{2}-\frac{2}{n}r(R+\triangle
\bar{f})]e^{-\bar{f}}dv\nonumber\\
&\geq&\int_{M}[\frac{2}{n}(R+\triangle\bar{f})^{2}-\frac{2}{n}r(R+\triangle
\bar{f})]e^{-\bar{f}}dv\nonumber\\
&\geq&\frac{2}{n}(\int_{M}(R+\triangle\bar{f})e^{-\bar{f}}dv)^{2}-
\frac{2r}{n}\int_{M}(R+\triangle\bar{f})e^{-\bar{f}}dv\nonumber\\
&=&\frac{2}{n}(\int_{M}(R+|\nabla\bar{f}|^{2})e^{-\bar{f}}dv)^{2}-
\frac{2r}{n}\int_{M}(R+|\nabla\bar{f}|^{2})e^{-\bar{f}}dv\nonumber\\
&=&\frac{2}{n}\lambda(t)(\lambda(t)-r).\nonumber
\end{eqnarray}
Now $\lambda(t)\leq r(t)$ and the monotonicity when
$\lambda(t)\leq0$ follows.

If $\frac{d}{dt}\lambda(t)=0$, then the equalities in the above
estimate hold. So we have $$R+\triangle\bar{f}=const.$$ and
$$Ric+\nabla^{2}\bar{f}=\frac{1}{n}(R+\bar{f})g.$$ On the other hand,
$e^{-\bar{f}/2}$ is the only eigenfunction of $-4\triangle+R$, so
$$2\triangle\bar{f}-|\nabla\bar{f}|^{2}+R=\lambda(t).$$ Thus
$\triangle\bar{f}-|\nabla\bar{f}|^{2}=const.,$ which equals zero
since
$\int_{M}(\triangle\bar{f}-|\nabla\bar{f}|^{2})e^{-\bar{f}}dv=0$.
By the maximal principle $f=const.$ over $M$. Hence $R=const.$ and
$Ric=\frac{R}{n}g$. This proves the lemma.
\end{proof}

\begin{lemma}[Perelman]\label{106}
Let $(M,g(t)),t\in[0,T),$ be a solution to the Ricci flow equation
(\ref{0}) or (\ref{101}). If $\lambda(g(t))>0$ for all $t$, then
$\nu(g(t))$ increases. Furthermore, the increasing is strict
unless the solution is a shrinking Ricci soliton.
\end{lemma}
\begin{proof}
We only need to prove the unnormalized case, since the $\nu$
functional is invariant up to rescalings and diffeomorphism
transformations, $i.e.$, $\nu(\alpha\phi^{*}g)=\nu(g)$ for any
constant $\alpha>0$ and diffeomorphism $\phi$ of $M$. Now for any
times $t_{1}$ and $t_{2}$ such that $0\leq t_{1}<t_{2}< T$, choose
$\tau_{0}>0$ and $f_{0}\in C^{\infty}$ satisfying
$\nu(g(t_{2}))=\mu(g(t_{2}),\tau_{0})=\mathcal{W}(g(t_{2}),f_{0},\tau_{0})$.
Solving the equation \cite[Equ. (3.3)]{P1}
\begin{equation}
\left\{ \begin{array}{ll}
\frac{\partial}{\partial t}g=-2Ric, \\
\frac{\partial}{\partial t}f=-\triangle f+|\nabla
f|^{2}-R+\frac{2n}{\tau},\\
\frac{\partial}{\partial t}\tau=-1,\\
f(t_{2})=f_{0},\tau(t_{2})=\tau_{0},
\end{array} \right.\nonumber
\end{equation}
and computing directly \cite{KL}, under this evolving equation, we
have
\begin{equation}\nonumber
\frac{d}{dt}\mathcal{W}=\int_{M}2\tau|Ric+\nabla^{2}f-\frac{1}{2\tau}g|^{2}(4\pi\tau)^{-n/2}e^{-f}dv\geq0.
\end{equation}
Hence by the definitions of $\mu$ and $\nu$, we have
\begin{eqnarray}
\nu(t_{1})&\leq&\mu(g(t_{1}),\tau_{0}+t_{2}-t_{1})\nonumber\\
&\leq&\mathcal{W}(g(t_{1}),f(t_{1}),\tau_{0}+t_{2}-t_{1})\nonumber\\
&\leq&\mathcal{W}(g(t_{2}),f(t_{2}),\tau_{0})=\nu(t_{2}).\nonumber
\end{eqnarray}
The equality doesn't hold unless the integrand equals zero,
$i.e.$, it is a shrinking Ricci soliton. On the other hand, $\nu$
remains obviously constant on a Ricci soliton by the invariance of
this functional under changes by rescalings and diffeomorphism
transformations.
\end{proof}

\vskip 10mm


\section{4-dimensional non-singular solutions}


From now on we are concerned with Ricci flow on $4$-manifolds. We
will continue to use the same notations and conventions in $\S 2$.
We assume all closed manifolds have constant volume $1$.

\begin{lemma} Let $(M,g(t)),t\in[0,\infty)$ be a  solution to
(\ref{0}) on a closed $4$-manifold $M$.  If $|R(g(t))|\leq C$ and
$ \breve{R}(g(t))\leq -c< 0$, where $C$ and $c$ are constants
independent of $t$, then
$$\int_{0}^{\infty}\int_{M}|Ric\textordmasculine (g(t))|^{2}dvdt<\infty.$$
\end{lemma}

\begin{proof} Note that, for $t\in [0, \infty)$,
$$|R(g(t))|< C, \ \ \ {\rm and} \ \ \ \breve{R}(g(t))\leq -c< 0,$$ where $C$ and $c$ are  constants  independent of $t$.  By
  Lemma 2.7, we have
$$\int_{0}^{\infty}\int_{M}|R-r|dv dt<\infty.$$ From the equation
$$\frac{\partial}{\partial t}R=\triangle R
+2|Ric\textordmasculine|^{2}+\frac{2}{4}R(R-r),$$ we obtain
\begin{eqnarray}
\int_{0}^{\infty}\int_{M} 2|Ric\textordmasculine|^{2}dv
dt&=&\int_{0}^{\infty}\int_{M}\frac{\partial}{\partial t}R dv dt-
\frac{1}{2}\int_{0}^{\infty}\int_{M} R(R-r)dv dt\nonumber\\
&=&\int_{0}^{\infty}\frac{\partial}{\partial t}r dt+
\frac{1}{2}\int_{0}^{\infty}\int_{M} R(R-r)dv dt\nonumber\\
&\leq&\lim_{t\longrightarrow \infty
}\sup|r(g(t))-r_{0}|+\frac{C}{2}\int_{0}^{\infty}\int_{M}|R-r|dv dt\nonumber\\
&\leq& 2C+\frac{C}{2}\int_{0}^{\infty}\int_{M}|R-r|dv
dt<\infty.\nonumber
\end{eqnarray}
\end{proof}

Observe that $\breve{R}(g(t))\leq \lambda _M=\overline{\lambda}
_M$. Theorem 1.4 follows immediately from

\begin{lemma}
Let $M$ be a closed oriented  4-manifold $M$ and let
$\{g(t)\},t\in[0,\infty)$, be a solution to (\ref{0}).  If
$|R(g(t))|<C$ and $\breve{R}(g(t))\leq
  -c< 0$,
where $C$ and $c$ are  constants  independent of $t$, then
$$2\chi(M)\geq 3|\tau(M)|.$$
Furthermore, if $\overline{\lambda}_{M}< 0$, then
$$2\chi(M)- 3|\tau(M)|\geq \frac{1}{96\pi^{2}}\overline{\lambda}_{M}^{2}.$$
and any non-singular solution $\{g(t)\},t\in[0,\infty)$ does not
collapse.
\end{lemma}

\begin{proof}
 From Lemma 3.1 we have
$$\int^{m+1}_{m}\int_{M} |Ric\textordmasculine(g(t))|^{2}dvdt \longrightarrow 0,$$ when
$m\longrightarrow \infty$.

By the Chern-Gauss-Bonnet formula and the Hirzebruch signature
theorem, for any metric $g$ on $M$,
$$\chi(M)=
\frac{1}{8\pi^{2}}\int_{M}(\frac{R(g)^{2}}{24}+|W^{+}(g)|^{2}+|W^{-}(g)|^{2}-\frac{1}{2}
|Ric\textordmasculine(g)|^{2})dv, \ \ \ \ { \rm and}$$
$$\tau(M)=\frac{1}{12\pi^{2}}\int_{M}(|W^{+}(g)|^{2}-|W^{-}(g)|^{2})dv,$$
where $W^{+}(g)$ and $W^{-}(g)$ are the self-dual and
anti-self-dual Weyl tensors respectively (cf. [1]).  Thus
\begin{eqnarray*} 2\chi(M)-3|\tau(M)|& \geq &
\liminf\limits_{m\longrightarrow \infty }
\frac{1}{4\pi^{2}}\int^{m+1}_{m}\int_{M}(\frac{1}{24}
R(g(t))^{2}-\frac{1}{2}|Ric\textordmasculine(g(t))|^{2})dvdt
\\ & = & \liminf\limits_{m\longrightarrow \infty }
\frac{1}{4\pi^{2}}\int^{m+1}_{m}\int_{M}\frac{1}{24}R(g(t))^{2}
dvdt \ge 0.\end{eqnarray*} This proves the first  inequality.

Observe that $\breve{R}(g(t))\leq \lambda _M=\overline{\lambda}
_M$. Note that
$$\int^{m+1}_{m}\int_{M}R(g(t))^{2}dvdt\geq
(\int^{m+1}_{m}r(g(t))dt)^{2}=(\int^{m+1}_{m}(r-\breve{R})dt+\int^{m+1}_{m}\breve{R}dt)^{2}.$$
From Lemma 2.7, we have
$$\int_{0}^{\infty}(r-\breve{R})dt<\infty.$$ Thus $$\lim\limits_{m\longrightarrow \infty
}\int^{m+1}_{m}(r-\breve{R})dt=0,$$ as $r-\breve{R}\geq 0$. By
taking $m\gg 1$ so that $|\int^{m+1}_{m}(r-\breve{R})dt|\ll 1$,
and
$$(\int^{m+1}_{m}(r-\breve{R})dt+\int^{m+1}_{m}\breve{R}dt)^{2}\geq
 (\int^{m+1}_{m}(r-\breve{R})dt+\overline{\lambda}_{M})^{2}.$$
 Thus we obtain  $$2\chi(M)-3|\tau(M)|\geq \liminf\limits_{m\longrightarrow \infty }
\frac{1}{96\pi^{2}}(\int^{m+1}_{m}(r-\breve{R})dt+\overline{\lambda}_{M})^{2}=
\frac{1}{96\pi^{2}}\overline{\lambda}_{M}^{2}.$$

The second inequality clearly implies that $\chi (M)>0$ whenever
$\overline{\lambda }_M<0$. On the other hand, by Cheeger-Gromov's
collapsing theorem  (c.f. [2] [3]) it holds that $\chi(M)=0$ if
$M$ collapse with bounded sectional curvature. This implies the
desired result.
\end{proof}

\begin{proof}[Proof of Theorem  \ref{002}] By Theorem 1.1 and the Hitchin-Thorpe inequality for
Einstein manifold it remains only to verify the inequality $$2\chi
(M)\ge 3|\tau (M)|$$ in the case (1.1.5). By $\S 2$ we know that
this only happens when $\breve{R}(g(t))\leq -c< 0$  and
$\lambda(g(t))\leq 0$, where $c$ is a constant independent of $t$.
By Lemma 3.2 the desired result follows.
\end{proof}

\begin{proof}[Proof of Corollary \ref{003}] If $M$ admits an $F$-structure of
positive rank, then $\chi(M)=0$ (c.f. [2] [3]). If $M$ admits a
shrinking soliton, by [8] the fundamental group $\pi _1(M)$ is
finite, and so $b_{1}(M)=b_{3}(M)=0$ by the Poincar\`e duality.
Hence $\chi(M)\geq 2$. By Theorem 1.2 the desired result follows.
\end{proof}

\begin{proof}[Proof of Corollary \ref{005}]By Theorem 1 in
   \cite{T}, the Spin$^{c}$ structure induced by a compatible almost complex
   structure on $(M, \omega)$ has Seiberg-Witten invariant equal to $\pm
   1$. Thus,  by Corollary 4.4 in \cite{Le},  we have
   $$\int_{M}R(g(t))^{2}dv\geq 32\pi^{2}(2\chi(M)+3\tau(M)).$$
   From the proof of Lemma 3.2, we know that \begin{eqnarray*}
2\chi(M)-3|\tau(M)| & \geq &  \liminf\limits_{m\longrightarrow
\infty }
\frac{1}{4\pi^{2}}\int^{m+1}_{m}\int_{M}\frac{1}{24}R(g(t))^{2}
dvdt
 \\ & \geq & \frac{1}{3}(2\chi(M)+3\tau(M)).\end{eqnarray*} Hence $$\chi(M) \geq 3\tau(M).$$

\end{proof}

\vskip 8mm

\section{Proof of Theorem 1.6}
\vskip 3mm

 Let $(M,g(t)),t\in[0,\infty)$ be a non-singular
solution to (\ref{0}) on a closed oriented  4-manifold $M$ with
$\overline{\lambda}_{M} < 0$. Since $|Rm(g(t))|<C$, there is a
constant $\varepsilon
>0$ depending only on $C$ such that, for any $t$,
$M_{t, \varepsilon}=\{x\in M: \text{Vol}(B_{x}(1), g(t))<
\varepsilon\}$ admits an
 F-structure of positive rank, and the Euler number  $\chi(M_{t, \varepsilon})=0$ (cf. \cite{CG1} \cite{CG2} \cite{A}).
By Theorem 1.4, $\chi(M)\geq
 \frac{1}{96\pi^{2}}\overline{\lambda}_{M}^{2}> 0$. Hence, for any $t$, there
 is an $x\in M$ such that $\text{Vol}(B_{x}(1), g(t))\geq \varepsilon$.

 By Lemma 2.7 and Lemma 3.1, we have  $$\int_{0}^{\infty}\int_{M} 2|Ric\textordmasculine|^{2}dv
dt< \infty, \ \ \ \int_{0}^{\infty}(r-\breve{R})dt<\infty \ \ \ \
{\rm and}$$ $$ \ \ \ \int_{0}^{\infty}\int_{M}|R-r|dvdt<\infty.$$
Thus we may choose a sequence of times $\{t_{k}\}$ so that
$t_{k}\longrightarrow
 \infty$ and  \begin{equation}\nu(k)=\int_{M} 2|Ric\textordmasculine(g(t_{k}))|_{k}^{2}dv_{k}\longrightarrow 0,
 \ \ \ \ |r(g(t_{k}))-\breve{R}(g(t_{k}))|\longrightarrow 0 \end{equation}  \begin{equation} \ \ \
\mu(k)=\int_{M}|R(g(t_{k}))-r(g(t_{k}))|dv_{k}\longrightarrow 0
\end{equation} as $k\longrightarrow \infty$.

We start the proof by proving three lemmas.

\begin{lemma} There is  a  sequence  points $\{x_{j,k}\in M\}$, $j=1, \cdots, m$,
satisfying that, for any $j$, $(M,g(t_{k}), x_{j, k})$
$C^{\infty}$-converges to a complete Einstein manifold
$(M_{j,\infty},g_{j,\infty},x_{j,\infty})$, i.e. there are
embeddings $F_{j,k, \rho }: B_{x_{j,\infty}}(\rho)\longrightarrow M$
such that, for any $\rho
> 0$, $F_{j,k,\rho}^{*}g(t_{k})$ $C^{\infty}$-converges to
$g_{j,\infty}$ on $B_{x_{j,\infty}}(\rho)\subset M_{j,\infty}$, and
$F_{j,k,\rho}( x_{j, \infty})=x_{j,k}$.  Furthermore,  all of the
manifolds $M_{j,\infty}$ are distinct, and, for any $\rho> 0$ and
$k\gg 1$, $\{F_{j,k,\rho}(B_{x_{j,\infty}}(\rho))\}$ are disjoint.
\end{lemma}

\begin{proof}  For each $g(t_k)$, choose a
 maximal number of disjoint unit balls $B_{x_{j,k}}(1)\subset (M, g(t_k))$
so that $\text{Vol}(B_{x_{j,k}}(1), g(t_{k}))\geq
 \varepsilon$. Since $\text{Vol}(M, g(t))\equiv 1$, there is a uniform bound on
 the number of the balls. Therefore, by passing to a subsequence if necessary, we may assume that
the maximal numbers of the disjoint unit balls for all  $g(t_{k})$
are the same, saying $\ell$.
  For every  $1\le j\le \ell$,  by Hamilton's compactness
theorem \cite{H1},  $(M,g(t_{k}+t), x_{j, k})$
$C^{\infty}$-converges to another normalized Ricci flow solution
$(M_{j,\infty},g_{j,\infty}(t),x_{j,\infty})$ after passing to a
subsequence.  By Proposition \ref{103},
$(M_{j,\infty},g_{j,\infty}(t),x_{j,\infty})$ is a complete
Einstein manifold with  negative scalar curvature for every $t$.
This proves that $(M,g(t_{k}), x_{j, k})$ $C^{\infty}$-converges
to a complete Einstein manifold
$(M_{j,\infty},g_{j,\infty},x_{j,\infty})$, where
$g_{j,\infty}\equiv g_{j,\infty}(0)$, i.e. there are
 embeddings  $F_{j,k, \rho }: B_{x_{j,\infty}}(\rho)\longrightarrow
M$ such that, for any $\rho > 0$, $F_{j,k,\rho}^{*}g(t_{k})$
$C^{\infty}$-converges to $g_{j,\infty}$ on
$B_{x_{j,\infty}}(\rho)\subset M_{j,\infty}$, and $F_{j,k,\rho}(
x_{j, \infty})=x_{j,k}$. Note that it is not necessary that all of
the manifolds $M_{j,\infty}$ are distinct. Let $\{M_{j,\infty}\}$,
$j=1, \cdots, m\le \ell $, be the resulting collection of distinct
manifolds. It is easy to see that, for any $\rho> 0$ and $k\gg 1$,
$\{F_{j,k,\rho}(B_{x_{j,\infty}}(\rho))\}$ are disjoint, and $\sum
\text{Vol}(M_{j,\infty}, g_{j,\infty})\leq 1$.
\end{proof}

\vskip 2mm

\noindent{\it Remark}: It is easy to see that two limit manifolds
$M_{j_1, \infty}$ and $M_{j_2, \infty}$ are the same if and only
if the distance $\text{dist}_{g(t_k)}(x_{j_1, k},x_{j_2, k})$ is
uniformly bounded above, independent of $k$.

\vskip 2mm

By the above, $(M_{j,\infty},g_{j,\infty},x_{j,\infty})$ satisfies
$|Rm(g_{j,\infty})|<C$ and $\text{Vol}(M_{j,\infty},
g_{j,\infty})\leq 1$. By \cite{CG3} there is a  good
   chopping
   of $M_{j,\infty}$, i.e. an exhaustion, $ \{U_{j,i}\}$, where every $ U_{j,i}$
is a compact $4$-submanifold with boundary $\partial U_{j,i}$,
$\bigcup_{i=1}^{\infty} U_{j,i}=M_{j,\infty}$ such that $$\cdots
\subset U_{j,i}\subset B_{x_{j,\infty}}(i) \subset
U_{j,i+1}\subset B_{x_{j,\infty}}(i+1) \subset  U_{j,i+2} \cdots
\subset
   M_{j,\infty}$$
satisfying that $ | \rm II \it (\partial U_{j,i}) |<\Lambda$ for
all $i$
   and $\rm Vol\it ( \partial U_{j,i}, g_{j,\infty}|_{\partial
U_{j,i}}) \longrightarrow 0$ as $i
   \longrightarrow \infty$, where $\rm II \it (\partial U_{j,i})$ is the
   second fundamental form of $\partial U_{j,i}$, and $\Lambda$ is
   a constant independent of $i$. By Lemma 4.1,
   we have   \begin{equation}\lim_{k\longrightarrow \infty}|\text{Vol}(U_{j,i}, F_{j,k,i+1}^{*}g(t_{k}))-
\text{Vol}(U_{j,i}, g_{j,\infty})|=0,\end{equation}
\begin{equation}\lim_{i\longrightarrow \infty}|\text{Vol}(U_{j,i},
g_{j,\infty})- \text{Vol}(M_{j,\infty}, g_{j,\infty})|=0, \ \ \ {\rm
and} \ \ \ \lim_{i\longrightarrow \infty}\text{Vol}(\partial
U_{j,i}, g_{j,\infty}|_{\partial U_{j,i}})=0.\end{equation}

\begin{lemma} Let $M_{k,i}=\coprod_{j=1}^{m}F_{j,k,i+1}(U_{j,i})\subset (M, g(t_k))$.  Then
$M_{k,i}$ is a $4$-submanifold of $(M, g(t_k))$  with boundary
   $\partial M_{k,i}\cong \coprod_{j=1}^{m}\partial U_{j,i}$, and every boundary component
is a graph 3-manifolds for $i\gg 1$.
\end{lemma}

\begin{proof} Since the
   second fundamental forms  $\rm II \it (\partial U_{j,i}) $ of $\partial
   U_{j,i}$ have a uniform bound, the sectional
   curvatures of $\partial U_{j,i}$ have a uniform bound, i.e. $|\rm Riem| _{\partial
   U_{j,i}}< \Lambda'$. Because  $\rm dim \partial U_{j,i}=3$ and $\lim_{i\longrightarrow \infty }\rm Vol\it ( \partial U_{j,i}, g_{j,\infty}|_{\partial
U_{j,i}})=0$, the components
   of $\partial U_{j,i}$ are graph 3-manifolds for $i\gg 1$ (cf. \cite{CG3}). By
   Lemma 4.1, $F_{j,k,i+1}(U_{j,i})$  are disjoint. Hence
   $M_{k,i}$ is a $4$-submanifold of $M$ with boundary $\partial M_{k,i} =\coprod_{j=1}^{m}\partial
   U_{j,i}$. The desired result follows.
\end{proof}

\begin{lemma} $$\lim_{k\longrightarrow
\infty}|\sum_{j=1}^{m}\rm {Vol}(U_{j,i}, F_{j,k,i+1}^{*}g(t_{k}))-
  \rm {Vol}(M, g(t_{k}))|\leq \overline{C}(\sum_{j=1}^{m} \rm Vol(\partial U_{j,i},
g_{j,\infty}|_{\partial U_{j,i}}))^{\frac{1}{2}},$$
$$\ \
\ {\rm and} \ \ \ \sum_{j=1}^{m} \rm Vol(M_{j,\infty},
 g_{j,\infty})= 1$$

where $\overline C$ is a constant independent of the indices.
\end{lemma}

\begin{proof} We first claim that there is an  $i_0>0$ such that, for any $i>i_0$, there
   is a $k_0$ satisfying that, for any
    $k> k_0$, $\rm Vol(B_{y}(1), g(t_{k}))\leq \varepsilon$ for all $y\in
   M- M_{k,i}$.

We may choose an $i_0\gg 1$ such that, for all
  $y_\infty \in \bigcup_{j=1}^{m}(M_{j,\infty}-
   B_{x_{j,\infty}}(i_0-2))$, $\rm Vol(B_{y_\infty}(1), g_{j,\infty})\leq
   \frac{1}{2}\varepsilon$. If the claim is false, for any fixed $i> i_0$,   there is a
  subsequence of times $\{t_{k_{s}}\}$, and a sequence of points
  $\{y_{k_{s}}\}$ such that  $y_{k_{s}}\in M- M_{k,i}$, and
\begin{equation} \rm Vol(B_{y_{k_{s}}}(1), g(t_{k_{s}}))>
\varepsilon
\end{equation}
Observe that the distance $\rm dist_{g(t_{k_{s}})}(y_{k_{s}},x_{j,
k_{s}})\longrightarrow \infty$ as $k_s\to \infty$ for all $1\le
j\le m$. Otherwise, assuming $\rm
dist_{g(t_{k_{s}})}(y_{k_{s}},x_{j, k_{s}})< {\rho }$ for some $j$
and $\rho >0 $, we get that $F_{j,k_{s},\rho
}^{-1}(y_{k_{s}})\longrightarrow y_{\infty}\in
B_{x_{j,\infty}}(\rho )- B_{x_{j,\infty}}(i-1)$, and so
\begin{equation}
\rm Vol(B_{y_{k_{s}}}(1), g(t_{k_{s}}))\longrightarrow \rm
Vol(B_{y_{\infty}}(1), g_{j,\infty})\le \frac 12 \varepsilon
\end{equation}
when $k_{s}\longrightarrow \infty$, since $F_{j,k_{s},
\rho}^{*}g(t_{k_{s}}) $ $C^{\infty}$-converges to $g_{j,\infty}$.
This contradicts to (11).

On the other hand, $(M,g(t_{k_{s}}), y_{ k_{s}})$
$C^{\infty}$-converges to a complete Einstein manifold
$(M_{\infty},g_{\infty},y_{\infty})$, and $M_{\infty}$ is distinct
from everyone of $M_{j,\infty}$ for $1\le j \le \ell$ (cf. the
remark after Lemma 4.1). This violates the choice of maximality of
$m$. The claim follows.

By \cite{CG3} and the above claim, for any $i>i_0$, there is a
$k_0$ such that, for any $k>k_0$,  $ M-M_{k,i}$ admits an
F-structure of positive rank, and $\chi(M- M_{k,i})=0$. By
Chern-Gauss-Bonnet theorem,  \begin{eqnarray*}0=\chi(M-
   M_{k,i})& =&
\frac{1}{8\pi^{2}}\int_{M-
M_{k,i}}(\frac{R(g(t_{k}))^{2}}{24}+|W^{\pm}(g(t_{k}))|_{k}^{2}-\frac{1}{2}
|Ric\textordmasculine(g(t_{k}))|_{k}^{2})dv_{k}\\
 & & +\int_{\partial
(M- M_{k,i})}P_{k}(\rm II\it), \end{eqnarray*} where $P_{k} (\rm
II\it) $ is a polynomial of the second
  fundamental form $\rm II\it(\partial
(M- M_{k,i}))$ of $\partial (M- M_{k,i})$ and its sectional
curvature. Hence
$$ \int_{M- M_{k,i}}\frac{R(g(t_{k}))^{2}}{24}dv_{k}\leq
\int_{M}\frac{1}{2} |Ric\textordmasculine(g(t_{k}))|_{k}^{2}dv_{k}
- 8\pi^{2}\int_{\partial (M- M_{k,i})}P_{k}(\rm II\it).$$
 Since $F_{j,k,i+1}^{*}g(t_{k})$ $C^{\infty}$-converges to
$g_{j,\infty}$ on $B_{x_{j,\infty}}(i+1)\subset M_{j,\infty}$, we
have $|P_{k}(\rm II\it)|< \overline{C}$, $k\gg 1$, where
$\overline{C}$ is a constant depending only on the bounds $\Lambda$
of  the
   second fundamental forms  $\rm II \it (\partial
U_{j,i})$ of $\partial U_{j,i}$, and
$$\upsilon(k): =|\rm Vol(\partial (M- M_{k,i}), g(t_{k})|_{\partial
(M- M_{k,i})})-\sum_{j=1}^{m} \rm Vol(\partial U_{j,i},
g_{j,\infty}|_{\partial U_{j,i}})|\longrightarrow 0$$ if
$k\longrightarrow \infty$. Hence
$$\int_{M- M_{k,i}}\frac{R(g(t_{k}))^{2}}{24}dv_{k}\leq
\frac{\nu(k)}{4}+\overline{C}(\sum_{j=1}^{m} \rm Vol(\partial
U_{j,i}, g_{j,\infty}|_{\partial U_{j,i}})+\upsilon(k))$$ where
$\nu (k)$ is defined in equation (7).

Clearly, \begin{eqnarray*}& & |r(g(t_{k}))|\rm Vol(M-
M_{k,i},g(t_{k}))-\int_{M- M_{k,i}}|R(g(t_{k}))|dv_{k}  \\
& \leq & \int_{M- M_{k,i}}|R(g(t_{k}))-r(g(t_{k}))|dv_{k}\leq
\mu(k), \ \ \
\end{eqnarray*}
where $\mu (k)$ is as in (8).

\begin{eqnarray*}\int_{M- M_{k,i}}|R(g(t_{k}))|dv_{k} & \leq &
(\int_{M- M_{k,i}}R(g(t_{k}))^{2}dv_{k})^{\frac{1}{2}}\rm Vol(M-
M_{k,i},g(t_{k}))^{\frac{1}{2}}\\ & \leq & (\int_{M-
M_{k,i}}R(g(t_{k}))^{2}dv_{k})^{\frac{1}{2}}.\end{eqnarray*} Thus
\begin{eqnarray*}\rm Vol(M-
M_{k,i},g(t_{k}))& \leq &
\frac{2}{|r_{\infty}|}(\mu(k)+(6\nu(k)\\
& & +\overline{C}(\sum_{j=1}^{m} \rm Vol(\partial U_{j,i},
g_{j,\infty}|_{\partial U_{j,i}})+\upsilon(k)))^{\frac{1}{2}}),
\end{eqnarray*}
where we used the fact that $r(g(t_{k}))\longrightarrow
r_{\infty}<0$ if $k\longrightarrow \infty$ (cf. the proof of
Proposition 2.6.) Therefore,

 \begin{eqnarray*}| \sum_{j=1}^{m}\rm Vol(U_{j,i},g_{j,\infty})-1|
& =& \lim_{k\longrightarrow \infty}|\rm Vol(M_{k,i}, g(t_{k}))-
  1| \\ & \leq &\overline{C}(\sum_{j=1}^{m} \rm Vol(\partial U_{j,i},
g_{j,\infty}|_{\partial U_{j,i}}))^{\frac{1}{2}}.\end{eqnarray*}
By letting $i\longrightarrow \infty$, we get the second equality
in the lemma.
\end{proof}

\begin{proof}[Proof of Theorem 1.6] By Lemma 4.3, for any $\delta >0$, we can choose $i\gg 1$ and $k\gg 1$ such
that $\rm Vol(M- M_{k,i},g(t_{k}))<\delta $.  Let $T=t_{k}$ and
$M^{\varepsilon}=M_{k,i}$. The desired results follows by Lemmas
4.1, 4.2 and 4.3.
\end{proof}

\vskip 8mm

\section{proof of Theorem 1.7}
\vskip 8mm

Let us first recall some facts about Seiberg-Witten equations,
which will
 be used to prove Theorem 1.7 (See \cite{Le} for details).
  Let $(M, g)$  be a   compact oriented  Riemannian  $4$-manifold
with a  $\rm Spin^{c}$ structure $\mathfrak{c}$. Let
$b^{+}_{2}(M)$ denote the dimension of the space of self-dual
harmonic $2$-forms in $M$. Let $S^{\pm}_{\mathfrak{c}}$ denote the
$\rm Spin^{c}$-bundles associated to $\mathfrak{c}$, and let $L$
be the determinant line bundle of $\mathfrak{c}$. There is a
well-defined Dirac operator
$$\mathcal{D}_{A}:
  \Gamma(S^{+}_{\mathfrak{c}})\longrightarrow
  \Gamma(S^{-}_{\mathfrak{c}})$$  Let $c:
  \wedge^{*}T^{*}M \longrightarrow {\rm End}(S^{+}_{\mathfrak{c}}\oplus
  S^{-}_{\mathfrak{c}})$ denote the Clifford multiplication on the
$\rm{Spin}^c$-bundles, and, for any $\phi\in \Gamma(S^{\pm})$, let
  $$q(\phi)=\overline{\phi}\otimes\phi-\frac{1}{2}|\phi|^{2}{\rm id}.$$
   The Seiberg-Witten equations read
   $$\begin{array}{ccc}\mathcal{D}_{A}\phi=0 \\
   c(F^{+}_{A})=q(\phi)
\end{array}   $$
where the unknowns are  a hermitian connection $A$ on $L$ and a
section $\phi\in \Gamma(S^{+}_{\mathfrak{c}})$, and $F^{+}_{A}$ is
the self-dual part of the curvature of $A$.  A resolution of the
Seiberg-Witten equations is called  {\it reducible} if $\phi\equiv
0$; otherwise, it is called {\it irreducible}.

 Let $(M,g(t)),t\in[0,\infty)$ be a non-singular
solution to (\ref{0}) on a closed oriented  4-manifold $M$ with
$\overline{\lambda}_{M} < 0$. We will continue to assume the
volume of $(M, g(t))$ is $1$. Let $\breve{R}(g(t))$ denote the
minimum of the scalar curvature of $g(t)$. Recall that
$\breve{R}(g(t))\le \overline{\lambda}_{M}$. Assume that $M$
  admits a symplectic structure $\omega$ satisfying that $b^{+}_{2}(M)
   >1$ and  $\chi(M)=3 \tau (M)$. Let $t_{k}$, $x_{j,k}$, $U_{j,i}$, $M_{j,\infty}$ and $F_{j,k,\rho}$
be the same as in Section  4.  By Theorem 1 in
   \cite{T}, the Spin$^{c}$ structure induced by a compatible almost complex
   structure $J$ on $(M, \omega)$ has Seiberg-Witten invariant equal to $\pm
   1$. Hence, for any $k$, there is an irreducible solution $(\phi_{k}, A_{k})$
to the Seiberg-Witten equations  (cf. \cite{Le}). For the sake of
simplicity we will use $|\cdot |_k$ to denote the norm with
respect to the metric $g(t_k)$.

   \begin{lemma} $$8\int_{M}|F^{+}_{A_{k}}|_{k}^{2}dv_{k} \geq
\int_{M}R(g(t_{k}))^{2}dv_{k}-\int_{M}48\pi^{2}
|Ric\textordmasculine(g(t_{k}))|_{k}^{2}dv_{k},$$
$$\ \ \ {\rm and} \ \ \ \lim_{k\longrightarrow \infty}\int_{M}|\nabla^{k} F^{+}_{A_{k}}|_{k}^{2}dv_{k}=0.$$
   \end{lemma}

   \begin{proof}The Bochner formula implies that
$$0=\frac{1}{2}\Delta_{k}
|\phi_{k}|_{k}^{2}+|\nabla^{A_{k}}\phi_{k}|_{k}^{2}+\frac{R(g(t_{k}))}{4}|\phi_{k}|_{k}^{2}+\frac{1}{4}|\phi_{k}|_{k}^{4},$$
$$4\int_{M}|\nabla^{A_{k}}\phi_{k}|_{k}^{2}dv_{k}=-\int_{M}(R(g(t_{k}))|\phi_{k}|_{k}^{2}+|\phi_{k}|_{k}^{4})dv_{k}.$$
From the estimate $|\phi_{k}|_{k}^{2}\leq - \breve{R}(g(t_{k}))$
(cf. \cite{Le}) and $-R(g(t_k))\le -\breve{R}(g(t_{k}))<0$ we get
that
$$4\int_{M}|\nabla^{A_{k}}\phi_{k}|_{k}^{2}dv_{k} \leq \int _M \breve{R}(g(t_{k}))^{2}dv_k-\int_{M}|\phi_{k}|_{k}^{4}dv_{k}.$$
By the second equation of the Seiberg-Witten equations and
$\chi(M)=3\tau(M)$,
\begin{eqnarray*}\int_{M}|\phi_{k}|_{k}^{4}dv_{k} & =&
8\int_{M}|F^{+}_{A_{k}}|_{k}^{2}dv_{k} \\ &\geq &
32\pi^{2}[c_{1}^{+}]^{2}[M] \geq 32\pi^{2}[c_{1}]^{2}[M] \\ &=&
32\pi^{2}(2\chi(M)+3\tau(M))\\ &=& 96 \pi^{2}(2\chi(M)-3\tau(M)) \\
&\geq & \int_{M}R(g(t_{k}))^{2}dv_{k}-\int_{M}48\pi^{2}
|Ric\textordmasculine(g(t_{k}))|_{k}^{2}dv_{k} ,
\end{eqnarray*}
where the last inequality follows by the Chern-Gauss-Bonnet
formula and Hirzebruch's signature formula (cf. section 4). Thus
\begin{eqnarray} 8\int_{M}|F^{+}_{A_{k}}|_{k}^{2}dv_{k} \geq
\int_{M}R(g(t_{k}))^{2}dv_{k}-\int_{M}48\pi^{2}
|Ric\textordmasculine(g(t_{k}))|_{k}^{2}dv_{k},\end{eqnarray}
where $c_{1}^{+}$ is the self-dual part of  the harmonic form
representing the first Chern class $c_{1}$ of $M$. Hence, by (7)
(8), we have
 \begin{eqnarray*}4\int_{M}|\nabla^{A_{k}}\phi_{k}|_{k}^{2}dv_{k} & \leq &
 \int_{M}(\breve{R}(g(t_{k}))^{2}-R(g(t_{k}))^{2})dv_{k}+\int_{M}48\pi^{2}
|Ric\textordmasculine(g(t_{k}))|_{k}^{2}dv_{k} \\ & \leq &
 C|\breve{R}(g(t_{k}))-r(g(t_{k}))|+\int_{M}48\pi^{2}
|Ric\textordmasculine(g(t_{k}))|_{k}^{2}dv_{k} \longrightarrow 0,
 \end{eqnarray*}$k\longrightarrow 0$,  where $C$ is a constant
 independent of $k$. By the second one of the Seiberg-Witten
 equations again (cf.
\cite{Le}), $$|\nabla^{k} F^{+}_{A_{k}}|_{k}^{2}\leq
 \frac{1}{2}|\phi_{k}|_{k}^{2}|\nabla^{A_{k}} \phi_{k}|_{k}^{2},$$ where $\nabla^{k}$ is the connection induced by
 Levi-civita connection.  Hence
 $$\int_{M}|\nabla^{k} F^{+}_{A_{k}}|_{k}^{2}dv_{k}\leq  \frac{1}{2}|\breve{R}(g(t_{k}))|\int_{M}|
\nabla^{A_{k}}\phi_{k}|_{k}^{2}dv_{k}
 \longrightarrow 0,$$ when $k\longrightarrow \infty$.
\end{proof}

Regard $F^{+}_{A_{k}}$ as self-dual 2-forms of  $g'(t_{k})$ on $
U_{j,i}\subset M_{j,\infty}$, where
 $g'(t_{k})=F_{j,k,i+1}^{*}g(t_{k})$. Since
 $|F^{+}_{A_{k}}|_{k}^{2}=\frac{1}{8}|\phi_{k}|_{k}^{4}\leq
 \frac{1}{8}\breve{R}(g(t_{k}))^{2}\leq C$,  where $C$ is a constant
 independent of $k$, $F^{+}_{A_{k}}\in L_{1}^{2}(g'(t_{k}))$, and $$\|F^{+}_{A_{k}}\|_{L_{1}^{2}(g'(t_{k}))}\leq C',$$
 where $C'$ is a constant independent of $k$. Note that $\|\cdot
\|_{L_{1}^{2}(g_{j,\infty})}\leq 2 \|\cdot
\|_{L_{1}^{2}(g'(t_{k}))}$ for $k\gg 1$ since $g'(t_{k})$
 $C^{\infty}$-converges to
$g_{j,\infty}$ on $ U_{j,i}$.   Thus, by passing to a subsequence,
$F^{+}_{A_{k}}$ $L_{1}^{2}$-converges to a 2-form $\Omega_{j}\in
L_{1}^{2}(g_{j,\infty})$, which is a  self-dual 2-form of
$g_{j,\infty}$.

  \begin{lemma} For any $j$, $\Omega_{j}$ is a smooth self-dual  2-form on
$U_{j,i}-
\partial U_{j,i}$ such that  $\nabla^{\infty} \Omega_{j}\equiv 0$,
and $| \Omega_{j}|_{\infty}\equiv {\rm cont.}\neq 0$, where
$\nabla^{\infty}$ is the connection induced by the
 Levi-civita connection of $g_{j,\infty}$.
    \end{lemma}

\begin{proof}
Note that
$$0 \leq \int_{U_{j,i}}|\nabla^{\infty} \Omega_{j}|_{\infty}^{2}dv_{\infty}=\lim_{k\longrightarrow \infty}
\int_{U_{j,i}}|\nabla^{\infty}
F^{+}_{A_{k}}|_{\infty}^{2}dv_{\infty} \leq \lim_{k\longrightarrow
\infty} \int_{M}|\nabla^{k} F^{+}_{A_{k}}|_{k}^{2}dv_{k}=0.$$  It
is easy to see that $\Omega_{j}$ is a weak solution of   the
elliptic equation $(d+d^{*})\Omega_{j}=0$ on $U_{j,i}$. By
elliptic equation theory, $\Omega_{j}$ is a smooth self-dual
2-form on $U_{j,i}-
\partial U_{j,i}$, and $\nabla^{\infty} \Omega_{j}\equiv 0$.

Now we claim that, for any $j$ and $i\gg 1$,
$\int_{U_{j,i}}|\Omega_{j}|_{\infty}^{2}dv_{\infty}\neq 0$. If it is
not true, there is a $j_{1}$ such that
$\int_{U_{j_{1},i}}|\Omega_{j_{1}}|_{\infty}^{2}dv_{\infty}\equiv
0$. Note that, by the results in Section 4,
\begin{eqnarray*}\int_{U_{j,i}}|\Omega_{j}|_{\infty}^{2}dv_{\infty}=\lim_{k\longrightarrow
\infty} \int_{U_{j,i}}| F^{+}_{A_{k}}|_{\infty}^{2}dv_{\infty}&=&
\lim_{k\longrightarrow \infty} \int_{U_{j,i}}|
F^{+}_{A_{k}}|_{k}^{2}dv_{k}\\
&\leq &\frac{1}{8}\lim_{k\longrightarrow \infty}
\breve{R}(g(t_{k}))^{2}\rm Vol(U_{j,i}, g'(t_{k}))\\ &=&
\frac{1}{8}r_{\infty}^{2}\rm Vol(U_{j,i}, g_{j,\infty})
\end{eqnarray*}
$$|\int_{M}(R(g(t_{k}))^{2}-r_{\infty}^{2})dv_{k}|\leq
 C\int_{M}(|R(g(t_{k})-r(g(t_{k}))|+|r_{\infty}-r(g(t_{k}))|)dv_{k}
 \longrightarrow 0,$$  and  \begin{eqnarray*}
\lim_{k\longrightarrow \infty}|\int_{M}|F^{+}_{A_{k}}|_{k}^{2}dv_{k}
-\sum_{j=1}^{m}
 \int_{U_{ji}}|F^{+}_{A_{k}}|_{k}^{2}dv_{k}|&=&\lim_{k\longrightarrow \infty}\int_{M- M_{k,i}}|F^{+}_{A_{k}}|_{k}^{2}dv_{k}
 \\& \leq & C \lim_{k\longrightarrow \infty}|\sum_{j=1}^{m}\rm Vol(U_{j,i},
 g'(t_{k}))\\ & &
  -
  \rm Vol(M, g(t_{k}))|\\ & \leq &\overline{C}(\sum_{j=1}^{m} \rm Vol(\partial U_{j,i},
g_{j,\infty}|_{\partial U_{j,i}}))^{\frac{1}{2}}, \end{eqnarray*} if
$k\longrightarrow \infty$, where $C$ and $\overline{C}$ are
constants independent of $k$.  By Lemma 5.1, we obtain
\begin{eqnarray*}r_{\infty}^{2}\sum_{j\neq j_{1}}\rm Vol(U_{j,i}, g_{j,\infty})&
\geq &
\sum_{j=1}^{m}\int_{U_{j,i}}8|\Omega_{j}|_{\infty}^{2}dv_{\infty}=\lim_{k\longrightarrow
\infty}\int_{M_{k,i}}8|F^{+}_{A_{k}}|_{k}^{2}dv_{k}\\ & \geq &
\lim_{k\longrightarrow
\infty}\int_{M}8|F^{+}_{A_{k}}|_{k}^{2}dv_{k}-\overline{C}(\sum_{j=1}^{m}
\rm Vol(\partial U_{j,i}, g_{j,\infty}|_{\partial
U_{j,i}}))^{\frac{1}{2}}
\\ & \geq& \lim_{k\longrightarrow \infty}\int_{M}(R(g(t_{k}))^{2}dv_{k}-\overline{C}(\sum_{j=1}^{m} \rm
Vol(\partial U_{j,i},
g_{j,\infty}|_{\partial U_{j,i}}))^{\frac{1}{2}}
\\ & =& r_{\infty}^{2}-\overline{C}(\sum_{j=1}^{m} \rm Vol(\partial U_{j,i},
g_{j,\infty}|_{\partial U_{j,i}}))^{\frac{1}{2}}.
\end{eqnarray*} Note that, for $i\gg 1$, $$1\gg \overline{C}(\sum_{j=1}^{m} \rm Vol(\partial U_{j,i},
g_{j,\infty}|_{\partial U_{j,i}}))^{\frac{1}{2}}\geq
r_{\infty}^{2}\rm Vol(U_{j_{1},i}, g_{j_{1},\infty}).$$  A
contradiction. Thus, for any $j$,
$\int_{U_{j,i}}|\Omega_{j}|_{\infty}^{2}dv_{\infty}\neq 0$. Thus
we obtain that $\nabla^{\infty} \Omega_{j}\equiv 0$, $|
\Omega_{j}|_{\infty}\equiv {\rm cont.}\neq 0$.
\end{proof}

\begin{proof}[Proof of Theorem 1.7]
  Since all $\Omega_{j}$, $1\le j\le m$, are self-dual 2-forms, and  $\nabla^{\infty} \Omega_{j}\equiv 0$, $|
\Omega_{j}|_{\infty}\equiv {\rm cont.}\neq 0$.  Hence, on any
$U_{j,i}$,  $g_{j, \infty}$ is a K\"{a}hler metric with K\"{a}hler
form $\sqrt{2}\frac{\Omega_{j}}{|\Omega_{j}|}$. By Lemma 4.1,
$g_{j, \infty}$ is a K\"{a}hler-Einstein  metric on
$M_{j,\infty}$.

Now, by the Chern-Gauss-Bonett theorem and the Hirzebruch theorem,
\begin{eqnarray*} 0=\chi(M)-3\tau(M) & \geq & \liminf_{k\longrightarrow
\infty}\frac{1}{2\pi^{2}}( \int_{M}|W^{-}(g(t_{k}))|_{k}^{2}dv_{k}
-\frac{1}{4} \int_{M}|Ric\textordmasculine(g(t_{k}))|_{k}^{2}dv_{k})\\
&= & \liminf_{k\longrightarrow \infty}\frac{1}{2\pi^{2}}
\int_{M}|W^{-}(g(t_{k}))|_{k}^{2}dv_{k} \geq 0,
\end{eqnarray*} where $W^{-}$ is the anti-self-dual Weyl tensor.
Thus $$\liminf_{k\longrightarrow \infty}
\int_{M}|W^{-}(g(t_{k}))|_{k}^{2}dv_{k} = 0.$$ Therefore, for any
$j$,
$$0\leq \int_{M_{j,\infty}}|W^{-}(g_{j,\infty})|_{\infty}^{2}dv_{\infty}\leq \liminf_{k\longrightarrow \infty}
\int_{M}|W^{-}(g(t_{k}))|_{k}^{2}dv_{k} = 0.$$ Hence
$g_{j,\infty}$ is a  K\"{a}hler-Einstein  metric with
$W^{-}(g_{j,\infty})\equiv 0$. This implies that $g_{j,\infty}$ is
a complex hyperbolic metric by the proof of Theorem 4.5 in
\cite{Le}. The desired result follows.
\end{proof}



\begin{thebibliography}{99}
\bibitem{A} M. T. Anderson, {\em Degeneration of metrics with bounded
curvature and applications to critical metrics of Riemannian
functionals}, Proceeding of Sympoia in Pure Mathematics, 54
(1993), 53-79.
\bibitem{B} A. L. Besse, {\em Einstein manifolds}, Ergebnisse der Math.
Springer-Verlag, Berlin-New York 1987.
\bibitem{CG1} J. Cheeger and M. Gromov, {\em Collapsing Riemannian Manifolds
 while keeping their curvature bounded I},  J.Diff.Geom. 23,
 (1986), 309-364.
\bibitem{CG2} J. Cheeger and M. Gromov, {\em Collapsing Riemannian Manifolds
 while keeping their curvature bounded II},  J.Diff.Geom. 32,
 (1990), 269-298.
 \bibitem{CG3} J.Cheeger, M.Gromov,  {\em On the characteristic numbers of
complete manifolds of bounded curvature and finite volume},
H.E.Rauch Memorial Volume I, Springer, Berlin, 1985, 115-154.
\bibitem{CGT} J. Cheeger, M. Gromov and M. Taylor, {\em Finite
propagation speed, kernal estimates for functions of the Laplace
operatro, and the geometry of complete Riemannian manifolds}, J.
Diff. Geom. 17 (1982), 15-53.
\bibitem{CZ} H.D. Cao and X.P. Zhu, {\em A complete proof of the
Poincar\'{e} and geometrization conjectures-application of the
Hamilton-Perelman theory of the Ricci flow}, Asian J. Math. 10
(2006) 165-492.
\bibitem{FZ} F. Fang and Y.G. Zhang., {\em Perelman's
$\lambda$-functional and the Seiberg-Witten equations}, preprint
2006
\bibitem{FG} M. Fern\'{a}ndez-L\'{o}pez and E. Garc\'{i}a-R\'{i}o,
{\em A remark on compact Ricci solitons}, preprint.

\bibitem{H2} R. Hamilton, {\em Three-manifolds with positive Ricci
curvature}, J. Diff. Geom. 17 (1982) 255-306.
\bibitem{H1} R. Hamilton, {\em A compactness property for solutions of
the Ricci flow}, Amer. J. Math. 117 (1995) 545-574.
\bibitem{H} R. Hamilton, {\em Non-singular solutions of the Ricci flow on
three-manifolds}, Comm. Anal. and Geom. 7 (1999) 695-729.
\bibitem{I} T. Ivey, {\em Ricci solitons on compact three-manifolds}, Diff. Geom.
Appl. 3 (1993) 301-307.
\bibitem{KL} B. Kleiner and J. Lott, {\em Notes on Perelman's Papers},
arXiv:math.DG/0605667.
\bibitem{Le}  C. LeBrun, {\em Four-Dimensional Einstein Manifolds and Beyond},
  in Surveys in Differential Geometry, vol VI: Essays on Einstein Manifolds,  247-285.
\bibitem{FIN} M. Feldman, T. Ilmanen and L. Ni, {\em Entropy and reduced distance for Ricci
expanders}, J. Geom. Anal. 15 (2005).
\bibitem{P1} G. Perelman, {\em The entropy formula
for the Ricci flow and its geometric applications},
arXiv:math.DG/0211159.
\bibitem{P2} G. Perelman, {\em Ricci flow with surgery on
three-manifolds}, arXiv:math.DG/0303109.
\bibitem{R} O. S. Rothaus, {\em Logarithmic Sobolev inequalities and
the spectrum of Schr$\ddot{\mbox{o}}$dinger operators}, J. Func.
Anal. 42 (1981) 110-120.
\bibitem{T} C. H. Taubes, {\em More constraints on symplectic forms from
 Seiberg-Witten invariants},  Math. Res. Lett. 2 (1995), 9-13.
\bibitem{SY}R.Schoen and S. T. Yau, {\em Lectures on differential
geometry}, in Conference Proceedings and Lecture Notes in Geometry
and Topology, 1, International Press Publications, 1994.
\end{thebibliography}
\end{document}